\documentclass[11pt]{article}
\usepackage{amsfonts}
\usepackage{bbm}
\usepackage{mathrsfs}
\usepackage{amscd}
\usepackage{color}
\usepackage{amsmath,amsfonts,amssymb,amscd}
\usepackage{indentfirst,graphics,epsfig,psfrag}
\input{epsf}
\usepackage{ifpdf}
\usepackage{enumerate}
\usepackage{appendix}
\usepackage{enumerate}

\usepackage{lineno}

\makeatletter \@addtoreset{figure}{section} \makeatother
\makeatletter
\long\def\@makecaption#1#2{%
   \vskip 10\p@
   \setbox\@tempboxa\hbox{{#1}\ \ #2}%
   \ifdim \wd\@tempboxa >\hsize

       {#1}\ \ #2\par
   \else
       \hbox to\hsize{\hfil\box\@tempboxa\hfil}%
   \fi}
\makeatother

\begin{document}
\title{\textbf{Normalized Laplacian spectrum
of some generalized subdivision-corona of two
regular graphs}}
\author{
\small  Qun Liu~$^{a,b}$, Zhongzhi Zhang~$^{a,c}$\\
\setcounter{footnote}{-1 }\thanks{Corresponding
author: Zhongzhi Zhang, E-mail address: liuqun09@yeah.net, zhangzz@fudan.edu.cn}\\
\small  a. School of Computer Science, Fudan University, Shanghai 200433, China\\
\small  b. School of Mathematics and Statistics, Hexi University, Gansu, Zhangye, 734000, P.R. China\\
\small  c. Shanghai Key Laboratory of Intelligent Information Processing, Fudan University, \\
\small  Shanghai 200433, China\\}
\date{}
\maketitle
\begin{abstract}
In this
paper, we determine the two normalized Laplacian spectrum of
generalized subdivision-vertex corona, subdivision-edge corona
for a connected regular graph with an arbitrary
regular graph
in terms of their normalized Laplacian eigenvalues. Moreover, applying
these results, we find some non-regular normaliaed Laplacian
cospectral graphs. These results generalize the existing results
in $\cite{AD}$.
\\[2mm]

{\bf Keywords:} Normalized Laplacian matrix; generalized
subdivision vertex corona; generalized subdivision edge coronae
\\[1mm]

{\bf MR(2010) Subject Classification:} 05C50/
{\bf CLC number:} O157.5
\end{abstract}

\section{ Introduction }
Spectral graph theory has an important role in determining
some structure and principal properties of a graph from its spectra.
There are several kinds of spectrum associated with a graph, for
example, adjacency spectrum, Laplacian spectrum, signless Laplacian
spectrum, normalized Laplacian spectrum etc. Normalized Laplacian
spectrum determines the bipartiteness from the largest eigenvalue
and the number of connected components from the second smallest
eigenvalue $\cite{FRK}$. Many researches have worked on the normalized Laplacian spectrum of
some graphs, the reader is referred to 
$(\cite{AP},\cite{SBJ},\cite{SB}-\cite{AD1})$.

Let $G$ be a simple and connected graph with vertex set
$V(G)$ and $E(G)$. Let $d_{i}$ be the degree of vertex $i$ in $G$ and $D_{G}=diag(d_{1},d_{2},\cdots d_{|V(G)|})$ the diagonal matrix with
all vertex degrees of $G$ as its diagonal entries. For a graph $G$,
let $A_{G}$ and $B_{G}$ denote the adjacency matrix and vertex-edge incidence matrix of $G$, respectively. The matrix $L_{G}=D_{G}-A_{G}$
is called the Laplacian matrix of $G$, where $D_{G}$ is the diagonal matrix of vertex degrees of $G$. F. Chung $\cite{FRK}$ introduced the normalized
Laplacian matrix of a graph $G$, denoted by $\mathcal{L}(G)$. It is
defined to be
\[
\begin{array}{crl}
\mathcal{L}(G)=I-D^{-\frac{1}{2}}AD^{-\frac{1}{2}}
=D^{-\frac{1}{2}}LD^{-\frac{1}{2}}.
\end{array}
\]
We denote the characteristic polynomial $det(xI-\mathcal{L}(G))$ of $\mathcal{L}(G)$ by $f_{G}(x)$. The
roots of $f_{G}(x)$ are known as the normalized Laplacian eigenvalues of $G$. The
multiset of the normalized Laplacian eigenvalues of $G$ is called the normalized Laplacian spectrum of $G$.
Chung $\cite{FRK}$ proved that all normalized Laplacian eigenvalues of a graph lie
in the interval $[0,2]$ and $0$ is always a normalized Laplacian
eigenvalue, that is $\lambda_{1}(G)=0$. She also determined normalized
Laplacian spectrum of different kinds of graph like complete graph, bipartite
graphs, hypercubes, etc. Two graphs $G$ and $H$ are called normalized
Laplacian cospectral or simply $\mathcal{L}$-cospectral if the spectrum of $\mathcal{L}(G)$
and $\mathcal{L}(H)$ are the same.

The subdivision graph $S(G)$ $\cite{CD}$
of a graph $G$ is the graph obtained by inserting a new
vertex into every edge of $G$. The set of such
new vertices is denoted by $I(G)$. In this paper, we are interested on finding normalized Laplacian spectrum of two generalized subdivision-corona of graphs, which are defined below.

{\bf Definition 1}\ Let $G_{1}$ and $H_{i}$ be two vertex-disjoint graphs
with number of vertices $n$ and $t_{i}$, and edges $m$ and $m_{i}$,
respectively. Then

(1)The generalized subdivision-vertex corona $\cite{LP}$ of $G$ and $H_{i}$ for $i=1,2,...,n$, denoted by $\Im(G)\odot\wedge_{i=1}^{n}H_{i}$,
is the graph obtained from $S(G)$ and $H_{i}$ by joining
the ith vertex of $V(G)$ to every vertex in $H_{i}$. The graph
$\Im(G)\odot\wedge_{i=1}^{n}H_{i}$ has $n+\sum_{i=1}^{n}t_{i}+m$ vertices and $2m+\sum_{i=1}^{n}m_{i}+\sum_{i=1}^{n}t_{i}$ edges.

(2) The generalized subdivision-edge corona of $G$ and
$H_{i}$ for $i=1,2,...,m$, denoted by $\Im(G)\ominus\wedge_{i=1}^{m}H_{i}$,
is the graph obtained from $S(G)$ and $H_{i}$ by joining
the ith vertex of $I(G)$ to every vertex in $H_{i}$. The graph
$\Im(G)\ominus\wedge_{i=1}^{m}H_{i}$ has $n+m+\sum_{i=1}^{m}t_{i}$
vertices and $2m+\sum_{i=1}^{m}t_{i}+\sum_{i=1}^{m}m_{i}$ edges.

Many researches have worked on the normalized Laplacian spectrum of
some graphs. In $\cite{CHY}$, Chen and Liao determined the normalized Laplacian spectrum of corona and edge corona of two graphs. In $\cite{AD,AD1}$, Das and $Panigrahi$ found the normalized Laplacian spectrum of some corona and subdivision-coronas of two regular graphs. Motivated by these works, here we determine
the normalized Laplacian spectrum of generalized subdivision-vertex corona
and subdivision-edge corona of two regular graph.

To prove our results, we need the following matrix products and
few results on them. For two matrices $A$ and $B$, of same size $m\times n$, the
Hadamard product $A\bullet B$ of $A$ and $B$ is a matrix of the
same size $m\times n$ with entries given by $(A\bullet B)_{ij}
=(A)_{ij}\cdot (B)_{ij}$(entrywise multiplication). Hadamard product is commutive, that is $A\bullet B=B\bullet A$. The line
graph $\cite{GC}$ of a graph $G$ is the graph $\emph{G}$, whose vertices are the edges of $G$ and two vertices of $\emph{l}(G)$ are adjacent if and only if they are incident on a common vertex in $G$. It is well known $\cite{CD}$ that $R^{T}(G)R(G)=A(\emph{l}(G))+2I_{m}$.



{\bf Lemma 1.1} $(\cite{CD})$ \ Let $M_{1}$, $M_{2}$, $M_{3}$
and $M_{4}$ be $p\times p$, $p\times q$, $q\times p$ and $q\times q$
matrices with $M_{1}$ and $M_{4}$ invertible. Then
\[
\begin{array}{crl}
 det\Bigg(
  \begin{array}{cccccccccccccccc}
    M_1 & M_2  \\
   M_3 & M_4  \\
 \end{array}
  \Bigg)&=&det(M_{4})det(M_{1}-M_{2}M_{4}^{-1}M_{3})\\
  &=&det(M_{1})det(M_{4}-M_{3}M_{1}^{-1}M_{2}),
\end{array}
\]
where $M_{1}-M_{2}M_{4}^{-1}M_{3}$ and $M_{4}-M_{3}M_{1}^{-1}M_{2}$
are called the Schur complements of $M_{4}$ and $M_{1}$, respectively.

{\bf Lemma 1.2} $(\cite{AD})$ \ For an $r$-regular graph $G$, the
eigenvalues of $A(\emph{l}(G))$ are the eigenvalues of $2(r-1)I_{n}-r\mathcal{L}(G)$ and $-2$ repeated $m-n$ times.

\section{Our results}
In this section, we focus on determing the
normalized Laplacian spectrum of generalized subdivision-vertex corona, subdivision-edge corona
whenever $G$ and $H_{i}$$(i=1,2,\ldots n)$ are respectively $r$-regular
and $r_{i}$$(i=1,2,\ldots,n)$-regular graph
in terms of Hadamard product of matrices.

{\bf Theorem 2.1} \ Let $G$ be an $r$-regular graph with $n$ vertices and $m$
edges, let $H_{i}$ be an $r_{i}$-regular graph with $t_{i}$ vertices for
$i=1,2,...,n$. Then we have the following:

(1)
\[
\begin{array}{lll}\mathcal{L}(\Im(G)\odot\wedge_{i=1}^{n}H_{i})&=&
\left(
  \scriptsize\begin{array}{ccccccc}
    I& -CR(G) &-M\\
    -CR^{T}(G)&I& 0\\
  -M^{T}& 0 & T^{'}\\
  \end{array}
\right),
\end{array}
\]
where $C$, $M$ and $T^{'}$ equal (2.1),(2.2)and (2.3), respectively.

(2)
\[
\begin{array}{lll}\mathcal{L}(\Im(G)\ominus\wedge_{i=1}^{m}H_{i})&=&
\left(
  \scriptsize\begin{array}{ccccccc}
    I& -CR(G) &0\\
    -CR^{T}(G)&I&-M\\
  0& -M^{T}& T^{'}\\
  \end{array}
\right),
\end{array}
\]
where $C$, $M$ and $T^{'}$ equal (2.4),(2.5)and (2.6), respectively.

{\bf Proof} \ (1) Let $R(G)$ be the incidence matrix
of $G$. The degree matrix and adjacency matrix of
$G^{*}=\Im(G)\odot\wedge_{i=1}^{n}H_{i}$ are
\begin{eqnarray*}
 D(G^{*})=\left(
  \begin{array}{cccccccccccccccc}
   P& 0 & 0\\
  0&2I_{m}&0\\
  0&0& T_{1}\\
 \end{array}
  \right),~~~~~~A(G^{*})=
  \left(
  \begin{array}{cccccccccccccccc}
   0& R(G)& Q\\
  R^{T}(G)&0&0\\
  Q^{T}&0&T_{2}\\
 \end{array}
  \right),
\end{eqnarray*}
where
\begin{eqnarray*}
 P=\left(
  \begin{array}{cccccccccccccccc}
   r+t_{1}& 0 &...&0\\
  0&r+t_{2}&...&0\\
  0&0& ...&0\\
   0&0&...&r+t_{n}\\
 \end{array}
  \right),~~~~~~Q=
  \left(
  \begin{array}{cccccccccccccccc}
   1^{T}_{t_{1}}& 0 &...&0\\
  0&1^{T}_{t_{2}}&...&0\\
  0&0&...&0\\
   0&0&...&1^{T}_{t_{n}}\\
 \end{array}
  \right),
\end{eqnarray*}

\begin{eqnarray*}
 T_{1}=\left(
  \begin{array}{cccccccccccccccc}
   (r_{1}+1)I_{t_{1}}& 0 &...&0\\
  0&(r_{2}+1)I_{t_{2}}&...&0\\
  0&0& ...&0\\
   0&0&...&(r_{n}+1)I_{t_{n}}\\
 \end{array}
  \right),~~~~~~T_{2}=
  \left(
  \begin{array}{cccccccccccccccc}
  A(H_{1})& 0 &...&0\\
  0&A(H_{2})&...&0\\
  0&0&...&0\\
   0&0&...&A(H_{n})\\
 \end{array}
  \right).
\end{eqnarray*}
Then
\[
\begin{array}{crl}
\mathcal{L}(G^{*})&=&I-D^{-\frac{1}{2}}(G)A(G)D^{-\frac{1}{2}}(G)\\
 &=&I-\left(
  \begin{array}{cccccccccccccccc}
   P^{-\frac{1}{2}}& 0& 0\\
  0& 2^{-\frac{1}{2}}I_{m}& 0 \\
  0& 0 & T_{1}^{-\frac{1}{2}}\\
 \end{array}
  \right)\left(
  \begin{array}{cccccccccccccccc}
   0& R(G)& Q\\
  R^{T}(G)& 0& 0 \\
  Q^{T}& 0 & T_{2}\\
 \end{array}
  \right)\left(
  \begin{array}{cccccccccccccccc}
   P^{-\frac{1}{2}}& 0& 0\\
  0& 2^{-\frac{1}{2}}I_{m}& 0 \\
  0& 0 & T_{1}^{-\frac{1}{2}}\\
 \end{array}
  \right)\\
  &=&\left(
  \begin{array}{cccccccccccccccc}
   I&  2^{-\frac{1}{2}}P^{-\frac{1}{2}}R(G)& P^{-\frac{1}{2}}QT_{1}^{-\frac{1}{2}}\\
  2^{-\frac{1}{2}}P^{-\frac{1}{2}}R^{T}(G)& I& 0 \\
  T_{1}^{-\frac{1}{2}}Q^{T}P^{-\frac{1}{2}}& 0 & I-T_{1}^{-\frac{1}{2}}T_{2}T_{1}^{-\frac{1}{2}}\\
 \end{array}
  \right).
\end{array}
\]
By computation, we have
\begin{eqnarray}C=2^{-\frac{1}{2}}P^{-\frac{1}{2}}
=diag\left(\frac{1}{\sqrt{2(r+t_{1})}},\frac{1}{\sqrt{2(r+t_{2})}},
\ldots,\frac{1}{\sqrt{2(r+t_{n})}}\right).
\end{eqnarray}

\begin{eqnarray}
 M=P^{-\frac{1}{2}}QT_{1}^{-\frac{1}{2}}
 =\left(
  \begin{array}{cccccccccccccccc}
 C^{T}_{t_{1}}& 0 & 0&...&0\\
  0&C^{T}_{t_{2}}&0&...&0\\
  0&0& ...&...&0\\
   0&0& 0&...&C^{T}_{t_{n}}\\
 \end{array}
  \right),
\end{eqnarray}
where $C^{T}_{t_{i}}=\frac{1}{\sqrt{(r+t_{i})(r_{i}+1)}}1^{T}_{t_{i}}$
$(i=1,2,\ldots,n)$.
\begin{eqnarray}
 T^{'}=I-T_{1}^{-\frac{1}{2}}T_{2}T_{1}^{-\frac{1}{2}}
 =\left(
  \begin{array}{cccccccccccccccc}
 L(H_{1})\bullet B(H_{1})&...&0\\
  0&L(H_{2})\bullet B(H_{2})&0\\
  0&...&0\\
   0&...&L(H_{n})\bullet B(H_{n})\\
 \end{array}
  \right).
\end{eqnarray}

(2) With the similar computation. The degree matrix and adjacency matrix of
$G^{\star}=\Im(G)\ominus\wedge_{i=1}^{m}H_{i}$ are
\begin{eqnarray*}
 D(G^{\star})=\left(
  \begin{array}{cccccccccccccccc}
   rI_{n}& 0 & 0\\
  0&P&0\\
  0&0& T_{1}\\
 \end{array}
  \right),~~~~~~A(G^{*})=
  \left(
  \begin{array}{cccccccccccccccc}
   0& R(G)& 0\\
  R^{T}(G)&0&Q\\
  0&Q^{T}&T_{2}\\
 \end{array}
  \right),
\end{eqnarray*}
where
\begin{eqnarray*}
 P=\left(
  \begin{array}{cccccccccccccccc}
   2+t_{1}& 0 &...&0\\
  0&2+t_{2}&...&0\\
  0&0& ...&0\\
   0&0&...&2+t_{n}\\
 \end{array}
  \right),~~~~~~Q=
  \left(
  \begin{array}{cccccccccccccccc}
   1^{T}_{t_{1}}& 0 & 0&...&0\\
  0&1^{T}_{t_{2}}&0&...&0\\
  0&0& ...&...&0\\
   0&0& 0&...&1^{T}_{t_{m}}\\
 \end{array}
  \right),
\end{eqnarray*}

\begin{eqnarray*}
 T_{1}=\left(
  \begin{array}{cccccccccccccccc}
   (r_{1}+1)I_{t_{1}}& 0 &...&0\\
  0&(r_{2}+1)I_{t_{2}}&...&0\\
  0&0& ...&0\\
   0&0&...&(r_{n}+1)I_{t_{n}}\\
 \end{array}
  \right),~~~~~~T_{2}=
  \left(
  \begin{array}{cccccccccccccccc}
  A(H_{1})& 0 &...&0\\
  0&A(H_{2})&...&0\\
  0&0&...&0\\
   0&0&...&A(H_{n})\\
 \end{array}
  \right).
\end{eqnarray*}
Then
\[
\begin{array}{crl}
\mathcal{L}(G^{\star})&=&I-D^{-\frac{1}{2}}(G^{\star})A(G^{\star})
D^{-\frac{1}{2}}(G^{\star})\\
 &=&I-\left(
  \begin{array}{cccccccccccccccc}
   r^{-\frac{1}{2}}I_{n}& 0 & 0\\
  0& P^{-\frac{1}{2}}& 0 \\
  0& 0 & T_{1}^{-\frac{1}{2}}\\
 \end{array}
  \right)\left(
  \begin{array}{cccccccccccccccc}
   0& R(G)& 0\\
  R^{T}(G)& 0& Q \\
  0& Q^{T}& T_{2}\\
 \end{array}
  \right)\left(
  \begin{array}{cccccccccccccccc}
  r^{-\frac{1}{2}}I_{n}& 0 & 0\\
  0& P^{-\frac{1}{2}}& 0 \\
  0& 0 & T_{1}^{-\frac{1}{2}}\\
 \end{array}
  \right)\\
  &=&\left(
  \begin{array}{cccccccccccccccc}
   I&  r^{-\frac{1}{2}}P^{-\frac{1}{2}}R(G)& 0\\
  r^{-\frac{1}{2}}P^{-\frac{1}{2}}R^{T}(G)& I& P^{-\frac{1}{2}}QT^{-\frac{1}{2}}_{1} \\
  0& T_{1}^{-\frac{1}{2}}Q^{T}P^{-\frac{1}{2}}& I-T_{1}^{-\frac{1}{2}}T_{2}T_{1}^{-\frac{1}{2}}\\
 \end{array}
  \right).
\end{array}
\]
By computation, we have
\begin{eqnarray}C=r^{-\frac{1}{2}}P^{-\frac{1}{2}}
=diag\left(\frac{1}{\sqrt{r(2+t_{1})}},\frac{1}{\sqrt{r(2+t_{2})}},
\ldots,\frac{1}{\sqrt{r(2+t_{n})}}\right).
\end{eqnarray}
\begin{eqnarray}
 M=P^{-\frac{1}{2}}QT_{1}^{-\frac{1}{2}}
 =\left(
  \begin{array}{cccccccccccccccc}
 C^{T}_{t_{1}}&  0&...&0\\
  0&C^{T}_{t_{2}}&...&0\\
  0&0& ...&0\\
   0& 0&...&C^{T}_{t_{n}}\\
 \end{array}
  \right),
\end{eqnarray}
where $C_{t_{i}}=\frac{1}{\sqrt{(2+t_{i})(r_{i}+1)}}$$(i=1,2,\ldots,n)$.

\begin{eqnarray}
 T^{'}=I-T_{1}^{-\frac{1}{2}}T_{2}T_{1}^{-\frac{1}{2}}
 =\left(
  \begin{array}{cccccccccccccccc}
\mathcal{L}(H_{1})\bullet B(H_{1})&...&0\\
  0&\mathcal{L}(H_{2})\bullet B(H_{2})&0\\
  0&...&0\\
   0&...&\mathcal{L}(H_{n})\bullet B(H_{n})\\
 \end{array}
  \right).
\end{eqnarray}

{\bf Notation}\ Let $G$ be a graph on $n$ vertices, $B$ and $C$
be matrices of size $n\times n$ and $n\times 1$, respectively.
For any parameter $x$, we have the notation: $\chi_{G}(B,C,x)=
C^{T}(xI_{n}-(L(G)\bullet B))^{-1}C$. We note that the notation
is similar to the notation coronal which was introduced by
Mcleman and McNicholas $\cite{MC}$.

{\bf Theorem 2.2} \ Let $G$ be an $r$-regular graph with $n$ vertices and $m$
edges, let $H_{i}$ be an $r_{i}$-regular graph with $t_{i}$ vertices for
$i=1,2,...,n$.
Then the normalized Laplacian spectrum of $\Im(G)\odot\wedge_{i=1}^{n}H_{i}$
consists of:

(1) The eigenvalue $\frac{1+r_{i}\delta_{ij}}{r_{i}+1}$
with multiplicity $n$, for every eigenvalue $\delta_{ij}$
$j=2,3,\ldots,t_{i}$ of $\mathcal{L}(H_{i})$,

(2) The eigenvalue $1$ with multiplicity $m-n$,

(3) The roots of the equation $2(r+t_{i}+rr_{i}+r_{i}t_{i})x^{3}$
$-2(2rr_{i}+2r_{i}t_{i}+3r+3t_{i})x^{2}+(2t_{i}r_{i}+4r+4t_{i}$
$+rr_{i}\mu_{i}+r\mu_{i})x-r\mu_{i}=0$ for each eigenvalue
$\mu_{i}$$(i=1,2\ldots, n)$ of $\mathcal{L}(G)$.

{\bf Proof} \ Let $R(G)$ be the incidence matrix
of $G$. With a suitable labeling for vertices,
the Normalized Laplacian characteristic polynomial of
$G^{*}=\Im(G)\odot\wedge_{i=1}^{n}H_{i}$ is
\[
\begin{array}{lll}f(G^{*},x)&=&
det\left(xI-\mathcal{L}(\Im(G)\odot\wedge_{i=1}^{n}H_{i})\right)\\
&=&det
\left(
  \scriptsize\begin{array}{ccccccc}
    (x-1)I_{n}& CR(G) &M\\
    CR^{T}(G)& (x-1)I_{m}& 0\\
  M^{T}& 0 & xI-T^{'}\\
  \end{array}
\right)\\
&=&det(xI-T^{'})det(S),
\end{array}
\]
where $C$, $M$ and $T^{'}$ equal (2.1),(2.2)and (2.3), respectively, and
\[
\begin{array}{crl}S&=&
\left(
  \begin{array}{cccccc}
    (x-1)I_{n}& CR(G) \\
 CR^{T}(G)& (x-1)I_{m}\\
  \end{array}
\right)-\left(
  \begin{array}{cc}
     M& \\
 0& \\
  \end{array}
\right)(xI-T^{'})^{-1}\left(
  \begin{array}{ccccc}
    M^{T}&0 \\
  \end{array}
\right)\\&=&\left(
  \begin{array}{ccccc}
    (x-1)I_{n}-M(xI-T^{'})^{-1}M^{T}& CR(G) &\\
 CR^{T}(G) & (x-1)I_{m}\\
  \end{array}
\right).\\
\end{array}
\]
Then
\begin{eqnarray*}
det(S)=(x-1)^{m}det\left((x-1)I_{n}-M(xI-T^{'})^{-1}M^{T}-
\frac{1}{x-1}C^{2}R(G)R(G)^{T}\right).
\end{eqnarray*}
~~~It is well known that $R(G)R^{T}(G)=A(G)+rI_{n}$ and $A(G)=r(I_{n}-\mathcal{L}(G))$,
so we get $R(G)R^{T}(G)=r(2I_{n}-\mathcal{L}(G))$.
Let $P_{i}=xI_{t_{i}}-L(H_{i})\bullet B(H_{i})$, then

$ M(xI-T^{'})^{-1}M^{T}$
\begin{eqnarray*}
 &=&\left(
  \begin{array}{cccccccccccccccc}
 C^{T}_{t_{1}}& 0 & 0&...&0\\
  0&C^{T}_{t_{2}}&0&...&0\\
  0&0& ...&...&0\\
   0&0& 0&...&C^{T}_{t_{n}}\\
 \end{array}
  \right)\left(
  \begin{array}{cccccccccccccccc}
 P^{-1}_{1}&...&0\\
  0&P^{-1}_{2}&0\\
  0&...&0\\
   0&...&P^{-1}_{n}\\
 \end{array}
  \right)\left(
  \begin{array}{cccccccccccccccc}
 C^{T}_{t_{1}}& 0 & 0&...&0\\
  0&C^{T}_{t_{2}}&0&...&0\\
  0&0& ...&...&0\\
   0&0& 0&...&C^{T}_{t_{n}}\\
 \end{array}
  \right)\\&=&\left(
  \begin{array}{cccccccccccccccc}
 \chi_{H_{1}}(B(H_{1}),C_{t_{1}},x)& 0 & 0&...&0\\
  0& \chi_{H_{2}}(B(H_{2}),C_{t_{2}},x)&0&...&0\\
  0&0& ...&...&0\\
   0&0& 0&...& \chi_{H_{n}}(B(H_{n}),C_{t_{n}},x)\\
 \end{array}
  \right).
\end{eqnarray*}
Therefore,
\[
\begin{array}{lll}det(S)&=&
(x-1)^{m}det\left((x-1)I_{_{n}}-\frac{C^{2}r(2I_{n}-\mathcal{L}(G))}{x-1}\right.
\\&&\left.-\left(
  \begin{array}{cccccccccccccccc}
 \chi_{H_{1}}(B(H_{1}),C_{t_{1}},x)& 0 & 0&...&0\\
  0& \chi_{H_{2}}(B(H_{2}),C_{t_{2}},x)&0&...&0\\
  0&0& ...&...&0\\
   0&0& 0&...& \chi_{H_{n}}(B(H_{n}),C_{t_{n}},x)\\
 \end{array}
  \right)\right)\\
 &=&(x-1)^{m-n}det\left((x-1)^{2}I_{_{n}}-C^{2}r(2I_{n}-\mathcal{L}(G))\right.
\\&&\left.-(x-1)\left(
  \begin{array}{cccccccccccccccc}
 \chi_{H_{1}}(B(H_{1}),C_{t_{1}},x)& 0 & 0&...&0\\
  0& \chi_{H_{2}}(B(H_{2}),C_{t_{2}},x)&0&...&0\\
  0&0& ...&...&0\\
   0&0& 0&...& \chi_{H_{n}}(B(H_{n}),C_{t_{n}},x)\\
 \end{array}
  \right)\right).
\end{array}
\]
~~~As $\mathcal{L}(H_{i})\bullet B(H_{i})=I_{t_{i}}-\frac{1}{r_{i}+1}A(H_{i})$,
we get $\mathcal{L}(H_{i})\bullet B(H_{i})=\frac{1}{r_{i}+1}(I_{t_{i}}
+r_{i}\mathcal{L}(H_{i}))$. So
\[
\begin{array}{lll}det(xI-T^{'})&=&
\left(
  \begin{array}{cccccccccccccccc}
 xI_{t_{1}}-\mathcal{L}(H_{1})\bullet B(H_{1})&...&0\\
  0&xI_{t_{2}}-\mathcal{L}(H_{2})\bullet B(H_{2})&0\\
  0&...&0\\
   0&...&xI_{t_{n}}-\mathcal{L}(H_{n})\bullet B(H_{n})\\
 \end{array}
  \right)\\
  &=&\prod_{i=1}^{n}\prod_{j=1}^{t_{i}}(x-\frac{1+r_{i}\delta_{ij}}{r_{i}+1}),
\end{array}
\]
where $\delta_{ij}$$(i=1,2,\ldots,n,j=1,2,\ldots,t_{i})$ is an eigenvalue of $\mathcal{L}(H_{i})$$(i=1,2,\ldots,n)$.

Since $H_{i}$$(i=1,2,\ldots,n)$ is regular, the sum of all entries on evergy row of its
normalized Laplacian matrix is zero. That means $\mathcal{L}(G_{i})C_{t_{i}}=(
1-\frac{r_{i}}{r_{i}})C_{t_{i}}=0C_{t_{i}}$. Then
$(\mathcal{L}(G_{i})\bullet B(G_{i}))C_{t_{i}}=(1-\frac{r_{i}}{r_{i}+1})C_{t_{_{i}}}=
\frac{1}{r_{i}+1}C_{t_{_{i}}}$ and $(xI_{t_{i}}-(\mathcal{L}(G_{i})\bullet B(G_{i}))
C_{t_{i}}=(x-\frac{1}{r_{i}+1})C_{t_{i}}$. Also, $C^{T}_{t_{i}}
C_{t_{i}}=\frac{t_{i}}{(r+t_{i})(r_{i}+1)}$.

Now $\chi_{H_{i}}(B(G_{i}),C_{t_{i}},x)=C^{T}_{t_{i}}(xI_{n_{i}}-(\mathcal{L}(G_{i})\bullet B(G_{i}))^{-1}C_{t_{i}}=\frac{C_{t_{i}}^{T}C_{t_{i}}}{x-\frac{1}{r_{i}+1}}=
\frac{t_{i}}{(r+t_{i})(r_{i}+1)(x-\frac{1}{r_{i}+1})}.$

Thus, if $\delta_{ij}$$(i=1,2,\ldots,n,j=1,2,\ldots,t_{i})$ is an eigenvalue of $\mathcal{L}(H_{i})$$(i=1,2,\ldots,n)$ and $\mu_{i}$
is an eigenvalue of $\mathcal{L}(G)$, then
\[
\begin{array}{lll}f_{G^{*}}(x)&=&
(x-1)^{m-n}\prod_{i=1}^{n}\prod_{j=1}^{t_{i}}(x-\frac{1+r_{i}\delta_{ij}}{r_{i}+1})\\&&
\prod_{i=1}^{n}\left[(x-1)(x-1-
\frac{t_{i}}{(r+t_{i})(r_{i}+1)(x-\frac{1}{r_{i}+1})})
+\frac{r(\mu_{i}-2)}{2(r+t_{i})}\right].
\end{array}
\]
(1) Since $H_{i}$$(i=1,2,\ldots,n)$ is connected, $0$ is a simple eigenvalue of $\mathcal{L}(H_{i})$. Also, since the only pole of $\chi_{G}(B(H_{i}),C_{t_{i}},x)$
is $x=\frac{1}{r_{i}+1},\frac{1+r_{i}\delta_{ij}}{r_{i}+1}$ an eigenvalue of $\mathcal{L}(\Im(G)\odot\wedge_{i=1}^{n}H_{i})$ with multiplicity $n$ for $j=2,3,\ldots,t_{i}$.

(2) Immediate from the characteristic polynomial.

(3) We get the remaining eigenvalues from the following equation:
$(x-1)(x-1$$-\frac{t_{j}}{(r+t_{j})(r_{j}+t_{j})(x-\frac{1}{r_{j}+1})})
+\frac{r(\mu_{i}-2)}{2(r+t_{i})}=0$, that is, $2(r+t_{i}+rr_{i}+r_{i}t_{i})x^{3}$
$-2(2rr_{i}+2r_{i}t_{i}+3r+3t_{i})x^{2}+(2t_{i}r_{i}+4r+4t_{i}$
$+rr_{i}\mu_{i}+r\mu_{i})x-r\mu_{i}=0$ for each eigenvalue
$\mu_{i}$$(i=1,2\ldots, n)$ of $\mathcal{L}(G)$.

{\bf Theorem 2.3} \ Let $G$ be an $r$-regular graph with $n$ vertices and $m$
edges, let $H_{i}$ be an $r_{i}$-regular graph with $t_{i}$ vertices for
$i=1,2,...,m$.
Then the normalized Laplacian spectrum of $\Im(G)\ominus\wedge_{i=1}^{m}H_{i}$
consists of:

(1) The eigenvalue $\frac{1+r_{i}\delta_{ij}}{r_{i}+1}$
with multiplicity $n$, for every eigenvalue $\delta_{ij}$
$(i=1,2,...,m, j=2,3,\ldots,t_{i})$ of $\mathcal{L}(H_{i})$,

(2) Two roots of the equation $(r_{i}t_{i}+2r_{i}+t_{i}+2)x^{3}-
(2r_{i}t_{i}+3t_{i}+4r_{i}+6)x^{2}+(r_{i}t_{i}
+r_{i}\mu_{i}+\mu_{i}+2t_{i}+4)x-\mu_{i}=0$
for each eigenvalue $\mu_{i}$$(i=1,2,...,n)$ of
$\mathcal{L}(G)$.

(3) Three roots of the equation $2(r+t_{i}+rr_{i}+r_{i}t_{i})x^{3}$
$(r_{i}t_{i}+2r_{i}+t_{i}+2)x^{3}-
(2r_{i}t_{i}+3t_{i}+4r_{i}+6)x^{2}+(r_{i}t_{i}
+r_{i}\mu_{i}+\mu_{i}+2t_{i}+4)x-\mu_{i}=0$ for each
eigenvalue $\mu_{i}$$(i=1,2,...,n)$ of $\mathcal{L}(G)$.

{\bf Proof} \ Let $R(G)$ be the incidence matrix
of $G$. With a suitable labeling for vertices,
the Normalized Laplacian characteristic polynomial of
$G^{*}=\Im(G)\ominus\wedge_{i=1}^{m}H_{i}$ is:
\[
\begin{array}{lll}f(\Im(G)\ominus\wedge_{i=1}^{m}H_{i},x)&=&
det(xI-\mathcal{L}(\Im(G)\ominus\wedge_{i=1}^{m}H_{i}))\\
&=&det
\left(
  \scriptsize\begin{array}{ccccccc}
    (x-1)I_{n}& CR(G) &0\\
    CR^{T}(G)&(x-1)I_{m}& M\\
  0& M^{T}& xI-T^{'}\\
  \end{array}
\right)\\
&=&det(xI-T^{'})detS,
\end{array}
\]
where
\[
\begin{array}{crl}S&=&
\left(
  \begin{array}{cccccc}
    (x-1)I_{n}& CR(G) \\
 CR^{T}(G)& (x-1)I_{m}\\
  \end{array}
\right)-\left(
  \begin{array}{cc}
     0& \\
 -M& \\
  \end{array}
\right)(xI-T^{'})^{-1}\left(
  \begin{array}{ccccc}
   0& -M^{T}\\
  \end{array}
\right)\\&=&\left(
  \begin{array}{ccccc}
    (x-1)I_{n}& CR(G) &\\
 CR^{T}(G) & (x-1)I_{m}-M(xI-T^{'})^{-1}M^{T}\\
  \end{array}
\right).\\
\end{array}
\]
Then
\begin{eqnarray*}
det(S)&=&(x-1)^{n}det\left((x-1)I_{m}-M(xI-T^{'})^{-1}M^{T}
-\frac{1}{x-1}C^{2}R(G)^{T}R(G)\right)\\
&=&(x-1)^{n}det\left((x-1)I_{m}-M(xI-T^{'})^{-1}M^{T}-\frac{C^{2}}{x-1}
(A(\emph{l}(G)+2I_{m})\right)\\
&=&(x-1)^{n}\left((x-1)I_{m}-M(xI-T^{'})^{-1}M^{T}
-\frac{C^{2}}{x-1}(-2+2)\right)^{m-n}
\\&&det\left((x-1-M(xI-T^{'})^{-1}M^{T}-\frac{C^{2}}{x-1}r(2I_{n}
-\mathcal{L}(G))\right)
(from~~Lemma\\
&=&\left((x-1)I_{m}-M(xI-T^{'})^{-1}M^{T}\right)^{m-n}
\\&&det\left((x-1)(x-1-M(xI-T^{'})^{-1}M^{T}-C^{2}r(2I_{n}
-\mathcal{L}(G))\right).
\end{eqnarray*}
Let $P_{i}=xI_{t_{i}}-L(H_{i})\bullet B(H_{i})$, then

$M(xI-T^{'})^{-1}M{T}$
\begin{eqnarray*}
 &=&\left(
  \begin{array}{cccccccccccccccc}
 C^{T}_{t_{1}}& 0 & 0&...&0\\
  0&C^{T}_{t_{2}}&0&...&0\\
  0&0& ...&...&0\\
   0&0& 0&...&C^{T}_{t_{n}}\\
 \end{array}
  \right)\left(
  \begin{array}{cccccccccccccccc}
 P^{-1}_{1}&...&0\\
  0&P^{-1}_{2}&0\\
  0&...&0\\
   0&...&P^{-1}_{n}\\
 \end{array}
  \right)\left(
  \begin{array}{cccccccccccccccc}
 C^{T}_{t_{1}}& 0 & 0&...&0\\
  0&C^{T}_{t_{2}}&0&...&0\\
  0&0& ...&...&0\\
   0&0& 0&...&C^{T}_{t_{n}}\\
 \end{array}
  \right)\\&=&\left(
  \begin{array}{cccccccccccccccc}
 \chi_{G_{1}}(B(H_{1}),C_{t_{1}},x)& 0 & 0&...&0\\
  0& \chi_{H_{2}}(B(H_{2}),C_{t_{2}},x)&0&...&0\\
  0&0& ...&...&0\\
   0&0& 0&...& \chi_{H_{n}}(B(G_{n}),C_{t_{n}},x)\\
 \end{array}
  \right).
\end{eqnarray*}

Hence,
\[
\begin{array}{lll}det(S)&=&
(x-1)^{n}det\left((x-1)I_{_{m}}-\frac{C^{2}R^{T}(G)R(G)}{x-1}\right.
\\&&\left.-\left(
  \begin{array}{cccccccccccccccc}
 \chi_{G_{1}}(B(H_{1}),C_{t_{1}},x)& 0 & 0&...&0\\
  0& \chi_{H_{2}}(B(H_{2}),C_{t_{2}},x)&0&...&0\\
  0&0& ...&...&0\\
   0&0& 0&...& \chi_{H_{n}}(B(G_{n}),C_{t_{n}},x)\\
 \end{array}
  \right)\right)\\
 &=&(x-1)^{n}det\left((x-1)^{2}I_{_{n}}-C^{2}(A(\emph{l}(G)+2I_{m}))\right.
\\&&\left.-\left(
  \begin{array}{cccccccccccccccc}
 \chi_{G_{1}}(B(H_{1}),C_{t_{1}},x)& 0 & 0&...&0\\
  0& \chi_{H_{2}}(B(H_{2}),C_{t_{2}},x)&0&...&0\\
  0&0& ...&...&0\\
   0&0& 0&...& \chi_{H_{n}}(B(G_{n}),C_{t_{n}},x)\\
 \end{array}
  \right)\right)
\end{array}
\]
As $\mathcal{L}(H_{i})\bullet B(H_{i})=I_{t_{i}}-\frac{1}{r_{i}+1}A(H_{i})$,
we get $\mathcal{L}(H_{i})\bullet B(H_{i})=\frac{1}{r_{i}+1}(I_{t_{i}}
+r_{i}\mathcal{L}(H_{i}))$.

So, if $\delta_{ij}$ is an eigenvalue of $\mathcal{L}(H_{i})$, then
\[
\begin{array}{lll}det(xI-T^{'})&=&
\left(
  \begin{array}{cccccccccccccccc}
 xI_{t_{1}}-\mathcal{L}(H_{1})\bullet B(H_{1})&...&0\\
  0&xI_{t_{2}}-\mathcal{L}(H_{2})\bullet B(H_{2})&0\\
  0&...&0\\
   0&...&xI_{t_{n}}-\mathcal{L}(H_{n})\bullet B(H_{n})\\
 \end{array}
  \right)\\
  &=&\prod_{i=1}^{m}\prod_{j=1}^{t_{i}}(x-\frac{1+r_{i}\delta_{ij}}{r_{i}+1}).
\end{array}
\]
~~~~Since $H_{i}$ is regular, the sum of all entries on every row of its
normalized Laplacian matrix is zero. That means $\mathcal{L}(H_{i})C_{t_{i}}=(
1-\frac{r_{i}}{r_{i}})C_{t_{i}}=0C_{t_{i}}$. Then
$(\mathcal{L}(H_{i})\bullet B(H_{i}))C_{t_{i}}=(1-\frac{r_{i}}{r_{i}+1})C_{t_{_{i}}}=
\frac{1}{r_{i}+1}C_{t_{_{i}}}$ and $(xI_{t_{i}}-(\mathcal{L}(H_{i})\bullet B(H_{i})))
C_{t_{i}}=(x-\frac{1}{r_{2}+1})C_{t_{i}}$. Also, $C^{T}_{t_{i}}
C_{t_{i}}=\frac{t_{i}}{(2+t_{i})(r_{i}+1)}$.

Now $\chi_{H_{i}}(B(G_{i}),C_{t_{i}},x)=C^{T}_{t_{i}}
\left(xI_{n_{i}}-(\mathcal{L}(G_{i})\bullet B(G_{i})\right)C_{t_{i}}=\frac{C_{t_{i}}^{T}C_{t_{i}}}{x-\frac{1}{r_{i}+1}}=
\frac{t_{i}}{(2+t_{i})(r_{i}+1)(x-\frac{1}{r_{i}+1})}.$

Thus, if $\delta_{ij}$ is an eigenvalue of $\mathcal{L}(H_{i})$ and $\mu_{i}$
is an eigenvalue of $\mathcal{L}(G)$, then
\[
\begin{array}{lll}f_{G^{*}}(x)&=&
(x-1-\frac{t_{i}}{(2+t_{i})(r_{i}+1)(x-
\frac{1}{r_{i}+1})})^{m-n}
\prod_{i=1}^{m}\prod_{j=1}^{t_{i}}
(x-\frac{1+r_{i}\delta_{ij}}{r_{i}+1})\\&&
\prod_{i=1}^{n}\left[(x-1)(x-1-
\frac{t_{i}}{(2+t_{i})(r_{i}+1)(x-\frac{1}{r_{i}+1})})
-\frac{2-\mu_{i}}{2+t_{i}}\right].
\end{array}
\]
(1) Since $H_{i}$$(i=1,2,\ldots,n)$ is connected, $0$ is a simple eigenvalue of $\mathcal{L}(H_{i})$. Also, since the only pole of $\chi_{G}(B(H_{i}),C_{t_{i}},x)$
is $x=\frac{1}{r_{i}+1},\frac{1+r_{i}\delta_{ij}}{r_{i}+1}$ an eigenvalue of $\Im(G)\ominus\wedge_{i=1}^{m}H_{i}$ with multiplicity $m$ for $j=2,3,\ldots,t_{i}$.

(2) We get $2(m-n)$ eigenvalues from the following equation:
$x-1-\frac{t_{i}}{(2+t_{i})(r_{i}+1)(x-
\frac{1}{r_{i}+1})}=0
$,
that is, $(2r_{i}+r_{i}t_{i}+t_{i}+2)x^{2}
-(2r_{i}+r_{i}t_{i}+2t_{i}+4)x+2=0$.

(3) The remaining eigenvalues can be obtained from the equation:
$(x-1)(x-1-
\frac{t_{i}}{(2+t_{i})(r_{i}+1)(x-\frac{1}{r_{i}+1})})
-\frac{2-\mu_{i}}{2+t_{i}}=0
$,
which is, $(r_{i}t_{i}+2r_{i}+t_{i}+2)x^{3}-
(2r_{i}t_{i}+3t_{i}+4r_{i}+6)x^{2}+(r_{i}t_{i}
+r_{i}\mu_{i}+\mu_{i}+2t_{i}+4)x-\mu_{i}=0$.

Now applying the results of this paper, we determine
some normalized Laplacian cospectral graphs. Since
for an $r$-regular graph $G$, we have $\mathcal{L}(G)=I_{n}-
\frac{1}{r}A(G)$, the lemma below is immediate.

{\bf Lemma 2.4}$(\cite{AD})$ \ Two regular graphs are $\mathcal{L}$-cospectral
if and only if they are cospectral.

By applying Theorem 2.2, we construct non-regular
$\mathcal{L}$-cospectral using generalized subdivision-
vertex corona and subdivision-edge corona.

{\bf Theorem 2.5} \ If $G_{1}$ and $H_{i}$$(i=1,2,...,n)$
are $\mathcal{L}$-cospectral regular graphs, and
$G_{2}$ and $H_{i}$$(i=1,2,...,n)$
are $\mathcal{L}$-cospectral regular graphs,
then $\Im(G)\odot\wedge_{i=1}^{n}H_{i}$ are $\mathcal{L}$-cospectral graphs.

{\bf Theorem 2.6} \ If $G_{1}$ and $H_{i}$$(i=1,2,...,m)$
are $\mathcal{L}$-cospectral regular graphs, and
$G_{2}$ and $H_{i}$$(i=1,2,...,m)$
are $\mathcal{L}$-cospectral regular graphs,
then $\Im(G)\ominus\wedge_{i=1}^{m}H_{i}$ are $\mathcal{L}$-cospectral graphs.

{\bf Remark 1}: In this paper, we determine the two normalized Laplacian spectrum of generalized subdivision-vertex corona, subdivision-edge corona for a connected regular graph with an arbitrary regular graph
in terms of their normalized Laplacian eigenvalues.These results generalize the existing results in $\cite{AD}$.

\vskip 0.1in
\noindent{\bf Acknowledgment:}
This work was supported by the National Natural Science Foundation of China (No. 11461020),


\begin{thebibliography}{99}

\bibitem{FRK} F.R.K. Chung, Spectral Graph Theory, in: CBMS. Reg. Conf. Ser. Math., vol. 92, AMS, providence, RI, 1997.
    
 \bibitem{ZZ} P.C. Xie, Z.Z. Zhang, Francesc Comellas, The normalized Laplacian spectrum of subdivisions of a graph, Appl. Math. Comput. 286 (2016): 250-256.
     
\bibitem{ZZ} P.C. Xie, Z.Z. Zhang, Francesc Comellas, On the spectrum of the normalized Laplacian of iterated triangulations of graphs, Applied Mathematics and Computation. 273 (2016): 1123-1129.
    
\bibitem{AP} A. Das, P. Panigrahi, Normalized Laplacian spectrum of some subdivision-joins and R-joins of two regular graphs, AKCE International Journal of Graphs and Combinatorics. (to appear) (2017).
    
 \bibitem{SBJ} S. Banerjee, J. Grout, On the spectrum of the normalized graph Laplacian, Linear Algebra Appl. 24 (2007): 47-56.
     
\bibitem{SB} S. Butler, J. Grout, A construction of cospectral graphs for the normalized Laplacian, Electronic journal of combinatorics. 18(1)(2011): paper 231.
    
\bibitem{HH} H. H. Li, J. S. Li, A note on the normalized Laplacian spectra, Taiwanese J. Math. 15(2011):129-139.
    
    

\bibitem{CD} Cvetkovi$\acute{c}$ D, Rowlinson P, Simi$\acute{c}$ S, An introduction to the theory of graph spectra. Cambridge: Cambridge University Press, 2009.

\bibitem{LP} P.L. Lu, Y.M. Wu, Laplacian and signless Laplacian characteristic polynomial of generalized subdivision corona vertex graph, Ars Combinatoria. 132 (2017): 357-369.


\bibitem{CHY}H.Y. Chen, L.W. L, The normalized Laplacian spectra of the corona and edge corona of two graphs, Linear and Multilinear Algebra. 65(3)(2017): 582-592.

\bibitem{AD}A. Das, P. Panigrahi, Normalized Laplacian spectrum of some subdivision-coronas of two regular graphs, Linear and Multilinear Algebra. 65(5)(2017): 962-972.

\bibitem{AD1}A. Das, P. Panigrahi, Normalized Laplacian spectrum of different type of corona of two regular graphs, Kragujevac Journal of Mathematics. 41(1)(2017): 57-69.



\bibitem{GC}Godsil C, Royle G. Algebraic graph theory.
New York(NY):Springer; 2001.


\bibitem{MC}McLeman C, McNocholas E, Spectra of coronae, Linear Algebra Appl. 435(2011): 998-1007.







\end{thebibliography}
\end{document}